\documentclass[12pt, a4paper, reqno]{amsart}

\usepackage[english]{babel}

\usepackage[T1]{fontenc}
\usepackage{lmodern}

\usepackage{amsmath}
\usepackage{amssymb}
\usepackage{amsthm}

\usepackage{hyperref}
\usepackage{url}

\usepackage{enumerate}
\usepackage[normalem]{ulem}

\renewcommand{\epsilon}{\varepsilon}
\renewcommand{\theta}[0]{\vartheta}
\renewcommand{\phi}[0]{\varphi}

\newcommand{\ppar}{\ \par}

\newcommand{\Span}[1]{\left\langle\, #1 \,\right\rangle}

\newcommand{\Set}[1]{\left\{ #1 \right\}}

\newcommand{\norm}[0]{\trianglelefteq}

\newcommand{\F}{\mathbb{F}}

\DeclareMathOperator{\GL}{GL}
\DeclareMathOperator{\SL}{SL}

\DeclareMathOperator{\Hom}{Hom}
\DeclareMathOperator{\End}{End}
\DeclareMathOperator{\Aut}{Aut}

\DeclareMathOperator{\Hol}{Hol}

\newtheorem{dummy}{Dummy}
\numberwithin{dummy}{section}
\numberwithin{figure}{section}

\newtheorem{theorem}[dummy]{Theorem}

\newtheorem{remark}[dummy]{Remark}

\theoremstyle{definition}

\theoremstyle{remark}

\makeatletter
\def\imod#1{\allowbreak\mkern10mu({\operator@font mod}\,\,#1)}
\makeatother

\numberwithin{equation}{section}

\usepackage{amsthm, mathtools, amssymb, amsfonts, latexsym, mathrsfs,
todonotes, eufrak, color}
\newcommand{\cC}{\mathcal{C}}

\begin{document}

\date{26 November 2019, 15:22 CET --- Version 1.02%
}

\title[Automorphism groups of groups of order $p^{2} q$]%
      {The automorphism groups of\\
        groups of order $p^{2} q$}
      
\author{E. Campedel}

\address[E.~Campedel]%
{Dipartimento di Matematica e Applicazioni\\
Edificio U5\\
Universit\`a degli Studi di Milano-Bicocca\\
via Roberto Cozzi, 55\\
20126 Milano}

\email{e.campedel1@campus.unimib.it}

\author{A. Caranti}

\address[A.~Caranti]%
 {Dipartimento di Matematica\\
  Universit\`a degli Studi di Trento\\
  via Sommarive 14\\
  I-38123 Trento\\
  Italy} 

\email{andrea.caranti@unitn.it} 

\urladdr{http://science.unitn.it/$\sim$caranti/}

\author{I. Del Corso}

\address[I.~Del Corso]%
        {Dipartimento di Matematica\\
          Universit\`a di Pisa\\
          Largo Bruno Pontecorvo, 5\\
          56127 Pisa\\
          Italy}
\email{delcorso@dm.unipi.it}

\urladdr{http://people.dm.unipi.it/delcorso/}



\begin{abstract}    
  We record for reference a detailed description  of the automorphism
  groups of the groups of order $p^{2}
  q$, where $p$ and $q$ are distinct primes.
\end{abstract}

\thanks{The second author is members of INdAM---GNSAGA. The second author
  gratefully acknowledges support from the Department of Mathematics of
  the University of Trento.}

\maketitle

\thispagestyle{empty}

\section{Introduction}

Let $p,  q$ be distinct  primes. O.~H\"older classified the  groups of
order $p^{2} q$ \cite{Hoe1}, and  also the groups of square-free order
\cite{Hoe2}.  H.~Dietrich and B.~Eick  \cite{DietEick} gave a detailed
description of the  structure of groups of cube-free  order, which was
complemented by  S.~Qiao and C.H.~Li \cite{QiaoLi}.   Groups of order
$p^{3}  q$ were  classified by  A.E.~Western \cite{West}  and R.~Laue
\cite{Laue}.   B.~Eick \cite{Eick4}  has given  an enumeration  of the
groups  whose order  factorises in  at  most $4$  primes. B.~Eick  and
T.~Moede \cite{EickMoede}  have enumerated  groups of order  $p^{n} q$,
for $n \le 5$.

For a forthcoming paper of ours, we need a detailed description of the
automorphism groups  of the groups of  order $p^{2} q$, where  $p$ and
$q$ are  distinct primes.  We  have recorded these data  for reference
here.

\section{The groups}

With current technology (i.e.~Sylow's theorem), describing the groups $G$
of order $p^{2} q$, where $p$ and $q$ are distinct primes, is an easy
exercise, which we now describe briefly, the basic point being that $G$ has
a normal Sylow subgroup.

If there are more than $1$, and thus exactly $q$, Sylow $p$-subgroups,
then $p  \mid q -  1$. 
\begin{itemize}
\item 
  If the  Sylow $p$-subgroups  intersect pairwise  trivially, counting
  $p$-elements show that then $G$ has exactly one Sylow $q$-subgroup.
\item 
  If there  are two distinct  Sylow $p$-subgroups $P_{1},  P_{2}$ that
  intersect  non-trivially in  a subgroup  $N$ of  order $p$,  then $N
  \norm G$, so  that $G$ has a  subgroup $R$ of order $p  q$, which is
  normal in $G$, as $p < q$. For the same reason, a Sylow $q$-subgroup
  of $R$ is normal  in $R$, and thus in $G$, so that  $G$ has a normal
  Sylow $q$-subgroup in this case as well.
\end{itemize}

We introduce some notation.
\begin{description}
  \item[$\cC_{n}$] denotes a cyclic group of order $n$.
  \item [$\rtimes$ and $\ltimes$] when they appear without subscripts,
    they denote the unique  (up to isomorphism) non-direct, semidirect
    product that is possible in the given situation.
  \item [$(\cC_{p} \times \cC_{p}) \rtimes_{S} \cC_{q}$] denotes a
    semidirect product where a generator of $\cC_{q}$ acts as a
    non-identity scalar matrix.
  \item[    $(\cC_{p} \times \cC_{p}) \rtimes_{D0} \cC_{q}$] denotes a
    semidirect product where a generator of $\cC_{q}$ acts as a
    diagonal, non-scalar matrix with no eigenvalue $1$, and
    determinant different from $1$.
  \item[    $(\cC_{p} \times \cC_{p}) \rtimes_{D1} \cC_{q}$] is the same as
    above, but the non-scalar matrix has
    determinant $1$ (and thus still no eigenvalue $1$).
  \item [$(\cC_{p} \times \cC_{p}) \rtimes_{C} \cC_{q}$] denotes a
    semidirect product where a generator of $\cC_{q}$ acts as a suitable
    power of a Singer cycle; note that the determinant of the matrix of a
    generator of $\cC_{q}$ acting on $\cC_{p} \times \cC_{p}$ is  $1$.
  \item [$\cC_{p^{2}} \ltimes_{1} \cC_{q}$] denotes a
    semidirect product with trivial centre.
  \item [$\cC_{p^{2}} \ltimes_{p} \cC_{q}$] denotes a
    semidirect product with  centre of order $p$.
\end{description}

Considering the  possible actions  on the
normal Sylow subgroup of  another Sylow subgroup, we obtain the following table. The automorphism
groups  are  determined in  Section~\ref{sec:auto},  on  the basis  of
results  of G.L.~Walls  \cite{Walls}, J.N.S.~Bidwell,  M.J.~Curran
and  D.J.~McCaughan  \cite{direct},  and  M.J.~Curran  \cite{Curran},
which we recall in Section~\ref{sec:Curran}.

\begin{center}
  \begin{tabular}{c|c|c|c|l}
    Type
    & Conditions & $G$ & $\Aut(G)$
    & Explanation
    \\\hline\hline
    1
    &

    &
    $\cC_{p^{2}} \times \cC_{q}$ & $\cC_{p} \times \cC_{p-1} \times
    \cC_{q-1}$
    &
    Cyclic groups
    \\\hline
    2
    &
    $p \mid q - 1$
    &
    $\cC_{p^{2}} \ltimes_{p} \cC_{q}$
    &
    $\cC_{p} \times \Hol(\cC_{q})$
    &
    Subs.~\ref{subs:p}
    \\\hline
    3
    &
    $p^{2} \mid q - 1$
    &
    $\cC_{p^{2}} \ltimes_{1} \cC_{q}$
    &
    $\Hol(\cC_{q})$
    &
    Thm~\ref{thm:Walls}
    \\\hline
    4
    &
    $q \mid p - 1$
    &
    $\cC_{p^{2}} \rtimes \cC_{q}$ 
    &
    $\Hol(\cC_{p^{2}})$
    &
    Thm~\ref{thm:Walls}
    \\\hline
    5
    &
    
    &
    $\cC_{p} \times \cC_{p} \times \cC_{q}$ 
    &
    $\GL(2, p) \times \cC_{q-1}$
    &
    Thm~\ref{thm:direct}
    \\\hline
    6
    &
    $q \mid p - 1$
    &
    $\cC_{p} \times (\cC_{p} \rtimes \cC_{q})$
    &
    $\cC_{p-1} \times \Hol(\cC_{p})$
    &
    Thms~\ref{thm:direct},~\ref{thm:Walls}
    \\\hline
    7 
    &
    $2 < q \mid p - 1$
    &
    $(\cC_{p} \times \cC_{p}) \rtimes_{S} \cC_{q}$
    &
    $\Hol(\cC_{p} \times \cC_{p})$
    &
    Subs.~\ref{subs:5-in-1},~\ref{subs:S}
    \\\hline
    8
    &
    $3 < q \mid p - 1$
    &
    $(\cC_{p} \times \cC_{p}) \rtimes_{D0} \cC_{q}$
    &
    $\Hol(\cC_{p}) \times \Hol(\cC_{p})$
    &
    Subs~\ref{subs:5-in-1},~\ref{subs:D}
    \\\hline
    9
    &
    $2 < q \mid p - 1$
    &
    $(\cC_{p} \times \cC_{p}) \rtimes_{D1} \cC_{q}$
    &
    $\cC_{2} \ltimes (\Hol(\cC_{p}) \times \Hol(\cC_{p}))$
    &
    Subs.~\ref{subs:5-in-1},~\ref{subs:D}
    \\\hline
    10
    &
    $2 < q \mid p + 1$
    &
    $(\cC_{p} \times \cC_{p}) \rtimes_{C} \cC_{q}$
    &
    $(\cC_{2} \ltimes \cC_{p^{2}-1}) \ltimes (\cC_{p} \times \cC_{p})$
    &
    Subs.~\ref{subs:5-in-1},~\ref{subs:C}
    \\\hline
    11
    &
    $p \mid q - 1$
    &
    $\cC_{p} \times (\cC_{p} \ltimes \cC_{q})$
    &
    $\Hol(\cC_{p}) \times \Hol(\cC_{q})$
    &
    Subs.~\ref{subs:last}
  \end{tabular}
\end{center}

\subsection{Isomorphism}

It is immediate to see that all types in this table
consist of exactly one isomorphism class of groups, with the exception of type
8. If $G$ is a group of this type, we can give it a canonical form by
choosing as generators first of all two eigenvectors with respect to distinct
eigenvalues in the normal, elementary abelian Sylow $p$-subgroup $V$. If
$\zeta$ is a fixed element of order $q$ in the multiplicative group of
the field with $p$ elements, we can then choose as a third generator
a suitable power $a$ of a $q$-element, so that it has eigenvalues
$\Set{\zeta, \zeta^{s}}$ on $V$. The parameter $s \notin \Set{0, 1,
  -1}$ determines $G$. If $t$ is the inverse of $s$ modulo $p$, then
$a^{t}$ has eigenvalues $ \Set{\zeta^{t}, \zeta}$ on $V$. It follows
that the parameters $s, t$ yield isomorphic groups, so that there are
$(q - 3)/2$ isomorphism classes of groups here.

\section{Automorphisms of (semi)direct products}
\label{sec:Curran}

We collect here the results we need of~\cite{direct, Walls, Curran}. We write
(auto)morphisms as exponents. 

\begin{theorem}[\protect{\cite[Theorem~3.2]{direct}}]
  \label{thm:direct}
  \ppar
  Let $G = H \times K$, where $H, K$ have no common direct factors.

  Then $\Aut(G)$ can be described in the natural way via the set of matrices
  \begin{align*}
    \left\{\ 
    \begin{bmatrix}
      a & c\\
      b & d
    \end{bmatrix}
    \right.
    :\ 
    &a \in \Aut(H), d \in \Aut(K),
    \\&b \in \Hom(K, Z(H)),\left. c \in \Hom(H, Z(K))
    \right\}.
  \end{align*}
\end{theorem}

\begin{theorem}\cite[Theorem 1]{Curran}
  \label{thm:semidirect}
  \ppar
  Let $G = H \rtimes K$ be a semidirect product.
  
  Then the subgroup of $\Aut(G)$  consisting of the automorphisms that
  leave $H$ invariant can be described in a natural way via the set of
  matrices
  \begin{align*}
    \left\{\ 
    \begin{bmatrix}
      a & 0\\
      b & d
    \end{bmatrix}
    \right.
    :\ &
    a \in \Aut(H), d \in \Aut(K),
    \\&
    (h^{k})^{a} = (h^{a})^{k^{d}} \text{, for $h \in H, k \in K$},
    \\&
    b : K \to H,
    (x y)^{b} = \left. x^{b} (y^{b})^{x^{d} } \text{, for $x, y \in K$}
    \ \right\}.
  \end{align*}
\end{theorem}

\begin{remark}
  \label{rem:a-and-d}
  The condition
  \begin{equation*}
    (h^{k})^{a} = (h^{a})^{k^{d}}
  \end{equation*}
  in Theorem~\ref{thm:semidirect}
  can be rewritten as
  \begin{equation*}
    \iota(k) ^{a} = a^{-1} \iota(k) a = \iota(k^{d}), 
  \end{equation*}
  where
  \begin{align*}
    \iota :\ &K \to \Aut(H)\\
    &k \mapsto (h \mapsto h^{k}).
  \end{align*}
  If $\Aut(H)$ is abelian, we get $\iota(k) = \iota(k^{d})$, that is,
  $[k, d] \in \cC_{K}(H)$. In particular, if $\cC_{K}(H) = 1$, then $d = 1$.
\end{remark}

\begin{theorem}[\protect{\cite[Theorem B]{Walls}, \cite[Example~1]{Curran}}]
  \label{thm:Walls}
  \ppar
  Let $G = \cC_{n} \rtimes \cC_{k}$, with $Z(G) = 1$. Write $H = \cC_{n}$, $K = \cC_{k}$.

  Then $H = G'$ is characteristic in $G$, and we have
  \begin{equation*}
    \Aut(G) \cong \Hol(\cC_{n}) = \cC_{n} \rtimes \Aut(\cC_{n}).
  \end{equation*}
\end{theorem}

\begin{remark}
  In    the    matrix     terms    of    Theorem~\ref{thm:semidirect},
  Theorem~\ref{thm:Walls} can be reformulated as
  \begin{align*}
    \Aut(G) =
    \left\{
      \begin{bmatrix}
        a & 0\\
        b & 1
      \end{bmatrix}
      :\right.\ &
      a \in \Aut(H), b : K \to H\\ 
      &
      (x y)^{b} =  x^{b} \left. (y^{b})^{x} \text{, for $x, y \in K$}
      \right\}.
  \end{align*}
  The $b$'s can be described in terms of the image $b_{0} \in H$ of a
  fixed generator of $K$: see Subsection~\ref{subs:5-in-1} for the details.
\end{remark}

\begin{theorem}[\protect{\cite[Theorem~3 and Example~1]{Curran}}]
  \label{thm:S}
  \ppar
  Let   $G   =   \cC_{n}    \rtimes   \cC_{k}$,   with   $Z(G)$   possibly
  non-trivial. Write $H = \cC_{n}$, $K = \cC_{k}$.

  Assume $H = G'$.

  Then 
  \begin{equation*}
  \Aut(G) \cong H \rtimes (\Aut(H) \times S),
  \end{equation*}
  where
  \begin{equation*}
  S
  = \Set{ d \in \Aut(K) : [k, d] = k^{-1} k^{d} \in \cC_{K}(H) \text{,
      for $k \in K$}}. 
  \end{equation*}
\end{theorem}

In matrix terms, Theorem~\ref{thm:S} states that
\begin{align*}
  \Aut(G)
  =
  \left\{\ 
  \begin{bmatrix}
    a & 0\\
    b & d
  \end{bmatrix}
  \right.
  :\ &
  a \in \Aut(H), d \in \Aut(K),
  \\&
    [k, d] = k^{-1} k^{d} \in \cC_{K}(H) \text{,
      for $k \in K$,}
    \\&\left. 
    b : K \to H, (x y)^{b} =  x^{b} (y^{b})^{x} \text{, for $x, y \in K$}
    \ \right\}
\end{align*}

\section{Automorphism groups}
\label{sec:auto}

We appeal to the results of
Section~\ref{sec:Curran}, whose
notation we employ. 

\subsection{Describing $b$}
\label{subs:5-in-1}
\label{subsub-b}

We begin by collecting some facts that hold true for
most cases.

Let us first consider the types 7, 8, 9, 10. Write
\begin{itemize}
\item 
  $\cC_{q} = \Span{z}$, 
\item
  $Z$ for the linear map $z$ induces
  on $H = \cC_{p} \times \cC_{p}$, and 
\item
  $Y$ for the linear map induced by $z^{d}$ on $H$. 
\end{itemize}
First note that for each $b_{0} \in  H$ there exists a unique function $b$
as in Theorem~\ref{thm:semidirect} such that $z^{b} = b_{0}$. In fact,
one has for $j = 1, \dots, q-1$
\begin{equation*}
  (z^{j})^{b} = b_{0}^{1 + Y + \dots + Y^{j-1}}.
\end{equation*}
In $\End(H)$ we have
\begin{equation*}
  0 = Y^{q} - 1
  =
  (Y - 1) (1 + Y + \dots + Y^{q-1}).
\end{equation*}
Now $Y  - 1$  invertible, as  $Y$ has no  eigenvalue $1$, so that $1 +
Y + \dots + Y^{q-1} = 0$. It follows that
\begin{equation*}
  (z^{q})^{b} 
= 
1
=
b_{0}^{1 + Y + \dots + Y^{q-1}},
\end{equation*}
is also satisfied.

A similar argument holds 
\begin{itemize}
\item
  for the types 4, 3, 2, 
\item
  for the subgroup $\cC_{p}
\rtimes \cC_{q}$ of type 6, and 
\item
  for the subgroup $\cC_{p} \ltimes \cC_{q}$ of type 11.
\end{itemize}

In these cases $Y$  is an automorphism of order coprime to
$r$ of a cyclic group $C$ of order  a power of a prime $r$, so that $Y
- 1$ is not nilpotent, and thus it is invertible, in $\End(C)$.

Note that
conjugating 
\begin{equation*}
  \begin{bmatrix}
    1 & 0\\
    b & 1\\
  \end{bmatrix}
\end{equation*}
by
\begin{equation*}
  \begin{bmatrix}
    a & 0\\
    0 & d\\
  \end{bmatrix}
\end{equation*}
we get
\begin{equation*}
  \begin{bmatrix}
    1 & 0\\
    d^{-1} b a & 1\\
  \end{bmatrix},
\end{equation*}
so that if $d = 1$ we have $z^{b a} = b_{0}^{a}$, and thus  the
group
\begin{equation*}
  \begin{bmatrix}
    a & 0\\
    b & 1\\
  \end{bmatrix}
\end{equation*}
is a split extension of $H$ by the group of the $a$'s.

\subsubsection{Between $d$ and $a$}
\label{subsub-d-and-a}

Suppose $d : z \mapsto z^{i}$, with $0 < i < q$ and $\gcd(i, q) = 1$.
For $h \in H$ we have
\begin{equation*}
  h^{a^{-1} Z a} = h^{Z^{i}},
\end{equation*}
and thus
\begin{equation}\label{eq:a}
  a^{-1} Z a = Z^{i}.
\end{equation}

\subsection{Type 7, $G = (\cC_{p} \times \cC_{p}) \rtimes_{S} \cC_{q}$}
\label{subs:S}

In this  case, since $Z$  is scalar,
we have $Z =  Z^{i}$, thus $q \mid i
- 1$, that  is, $i  = 1$ and  $d$ is
trivial. Since  $a$ is arbitrary, we
obtain as the automorphism group the
holomorph  of $\cC_{p}  \times \cC_{p}$,
that   this  the  affine   group  in
dimension $2$ over $\F_{p}$.

\subsection{Type 8 and 9, $G = (\cC_{p} \times \cC_{p}) \rtimes_{D0} \cC_{q}$ or $(\cC_{p} \times \cC_{p}) \rtimes_{D1} \cC_{q}$}
\label{subs:D}

In this case
\begin{equation*}
  Z =
  \begin{bmatrix}
    \lambda & 0\\
    0 & \mu\\
  \end{bmatrix},
\end{equation*}
with $\lambda \ne \mu$, $\lambda, \mu \ne 1$. Then~\eqref{eq:a} yields
$\Set{\lambda,  \mu}  = \Set{\lambda^{i},  \mu^{i}}$.   If $\lambda  =
\lambda^{i}$ and $\mu  = \mu^{i}$, we obtain that $q \mid  i - 1$, and
thus    $i    =    1$    and   $d    =    1$.    From~\eqref{eq:a}~and
Subsection~\ref{subsub-b}, we obtain  that $a$ centralizes $Z$, and
that the automorphism group contains
\begin{equation}\label{eq:group-D0}
  (\cC_{p-1} \times \cC_{p-1}) \ltimes (\cC_{p} \times \cC_{p})
  =
  \Hol(\cC_{p}) \times \Hol(\cC_{p}),
\end{equation}
with $\cC_{p-1} \times \cC_{p-1}$ acting by diagonal matrices on $\cC_{p}
\times \cC_{p}$, a 
typical element being
\begin{equation}\label{eq:elt-of-D0}
  \begin{bmatrix}
    T & 0\\
    b & 1\\
  \end{bmatrix}
\end{equation}
with $T$ diagonal.

If $\lambda =  \mu^{i}$ and $\mu = \lambda^{i}$, then $\lambda =
\lambda^{i^{2}}$, so that $q \mid (i - 1) (i + 1)$. When $q \mid i -
1$ we get again $d = 1$, whereas when $q \mid i + 1$ we get $z^{d} =
z^{-1}$ and $\lambda = \mu^{-1}$. Thus this case only occurs when
$\det(Z) = 1$, that is, when $G$ is of type 9. The inversion $d$ can then be paired with 
\begin{equation*}
  S =
  \begin{bmatrix}
    0 & 1\\
    1 & 0\\
  \end{bmatrix},
\end{equation*}
to get $S^{-1} Z S = Z^{-1}$. In this case the automorphism group is
the extension of the group~\eqref{eq:group-D0} by the involution
\begin{equation}\label{eq:inv-for-D0}
  \begin{bmatrix}
   S & 0\\
   0 & d\\
  \end{bmatrix},
\end{equation}
which acts on~\eqref{eq:elt-of-D0} as
\begin{equation*}
  \begin{bmatrix}
    S & 0\\
    0 & d\\
  \end{bmatrix}^{-1}
  \cdot 
  \begin{bmatrix}
    T & 0\\
    b & 1\\
  \end{bmatrix}
  \cdot
    \begin{bmatrix}
    S & 0\\
    0 & d\\
  \end{bmatrix}
  =
  \begin{bmatrix}
    S T S & 0\\
    d^{-1} b S & 1\\
  \end{bmatrix};
\end{equation*}
Now 
\begin{equation*}
  z^{d^{-1} b S}
  =
  (z^{-1})^{b S}
  =
  (z^{q-1})^{b S}
  =
  (b_{0}^{1 + Y + \dots + Y^{q-2}})^{S}
  =
  (b_{0}^{- Y^{q-1}})^{S}
  =
  b_{0}^{- Y^{-1} S},
\end{equation*}
where 
\begin{equation*}
  - Y^{-1} S
  =
  \begin{bmatrix}
    0 & - \lambda^{-1}\\
    - \lambda & 0\\
  \end{bmatrix}
\end{equation*}
is an involution, that acts by exchanging the
two copies of $\Hol(\cC_{p})$.

\subsection{Type 10, $G = (\cC_{p} \times \cC_{p}) \rtimes_{C} \cC_{q}$}
\label{subs:C}

Note first that the order $q \ne 2$ of $Z$ divides $p + 1$, so it does
not divide $p - 1$. It follows that $Z \in \SL(2, p)$, that is, $\det(Z) = 1$.

If $\lambda, \mu = \lambda^{-1}$ are the (distinct)
eigenvalues of $Z$ in the field
$\F_{p^{2}}$, then $\Set{\lambda, \mu} = \Set{\lambda^{i}, \mu^{i}}$.
If $\lambda = \lambda^{i}$ and $\mu = \mu^{i}$, we get once more $d =
1$. Thus in this case $a$ 
lies in the centralizer of $Z$ in $\Aut(H) = \GL(2, p)$ 
\begin{equation*}
  \cC_{\Aut(H)}(Z) 
  = 
  \Set{ u + v Z \ne 0 : u, v \in \F_{p} },
\end{equation*}
which is cyclic, of order $p^{2} - 1$.

If $\lambda =  \mu^{i}$ and $\mu = \lambda^{i}$, then $\lambda =
\lambda^{i^{2}}$, so that $q \mid (i - 1) (i + 1)$. When $q \mid i -
1$ we get again $d = 1$, whereas when $q \mid i + 1$ we get $z^{d} =
z^{-1} = z^{p}$, as $p \equiv -1 \pmod{q}$.

In an appropriate basis of $H$ we have
\begin{equation*}
  Z =
  \begin{bmatrix}
    0 & 1\\
    -1 & t\\
  \end{bmatrix},
\end{equation*}
where $t = \lambda + \lambda^{-1} = \lambda + \lambda^{p}$.

The equation $\iota(z)^{a} = \iota(z^{d})$ of Remark~\ref{rem:a-and-d}
has now become
\begin{equation*}
  Z^{a} = Z^{-1}.
\end{equation*}
One solution $a$ for this is
\begin{equation*}
  S =
  \begin{bmatrix}
    0 & 1\\
    1 & 0\\
  \end{bmatrix},
\end{equation*}
as
\begin{equation*}
  S^{-1} Z S
  =
  \begin{bmatrix}
    t & - 1\\
    1 & 0\\
  \end{bmatrix}
  =
  Z^{-1}.
\end{equation*}
All other solutions are  obtained in the form $a = X  C$, where $X \in
\cC_{\Aut(H)}(Z)$. So  in this  case $\Aut(G)$ is  the extension  of the
subgroup determined  by $a  = d  = 1$, which  is isomorphic  to $\cC_{p}
\times \cC_{p}$, acted  upon by the subgroup determined by  $b = 0$. The
latter subgroup has a normal subgroup
\begin{equation*}
  \mathcal{C}
  =
  \begin{bmatrix}
    \cC_{\Aut(H)}(Z) & 0\\
    0 & 1\\ 
  \end{bmatrix},
\end{equation*}
which is cyclic, of order $p^{2} - 1$
extended by the involution
\begin{equation*}
  D
  =
  \begin{bmatrix}
   S & 0\\
   0 & d\\
  \end{bmatrix},
\end{equation*}
where $z^{d} = z^{-1} = z^{p}$. Now $\cC_{\Aut(H)}(Z)$ is the multiplicative
group of the field with $p^{2}$ elements; conjugation by $S$ induces
an automorphism group of order $2$, which is then the Frobenius
map. Thus $\cC_{\mathcal{C}}(D) =
\cC_{\Aut(H)}(S)$ has order $p-1$, 
and $g^{S} = g^{p}$ for $g \in \mathcal{C}$. 

\subsection{Type 2, $G = \cC_{p^{2}} \ltimes_{p} \cC_{q}$}
\label{subs:p}

Here $\cC_{p^{2}} = \Span{x}$ induces on $\cC_{q}$ a group of
automorphisms of order $p$, and thus the centraliser $\cC_{\Span{x}}(\cC_{q}) =
\Span{x^{p}}$ has order $p$.

As per Remark~\ref{rem:a-and-d}, here
\begin{equation*}
  S = \Set{ 
    d \in \Aut(\cC_{p^{2}}) : 
    [x, d] \in \cC_{\Span{x}}(\cC_{q}) =
    \Span{x^{p}} }.
\end{equation*}
is a group of order $p$, generated by the automorphism $x \mapsto x^{1
  + p}$. According to Theorem~\ref{thm:S}, we
get that the automorphism group is isomorphic to 
\begin{equation*}
  \cC_{q} \rtimes (\cC_{q-1} \times \cC_{p})
  \cong
  \cC_{p} \times \Hol(\cC_{q}).
\end{equation*}

\subsection{Type 11, $G = \cC_{p} \times (\cC_{p} \ltimes \cC_{q})$}
\label{subs:last}

According to Theorem~\ref{thm:direct}~and Theorem~\ref{thm:Walls}, we
have that the automorphism group has the form
\begin{equation*}
  \begin{bmatrix}
    \cC_{p-1} & 0\\
    \cC_{p} & \Hol(\cC_{q})\\
  \end{bmatrix}.
\end{equation*}
To see the structure, let us consider the conjugate
\begin{equation*}
  \begin{bmatrix}
    a & 0\\
    0 & d\\
  \end{bmatrix}^{-1}
  \cdot
  \begin{bmatrix}
    1 & 0\\
    b & 1\\
  \end{bmatrix}
  \cdot
  \begin{bmatrix}
    a & 0\\
    0 & d\\
  \end{bmatrix},
\end{equation*}
where $a \in \Aut(\cC_{p})$, $d \in \Aut(\cC_{p} \ltimes \cC_{q}) \cong
\Hol(\cC_{q})$, and $b \in \Hom(\cC_{p} \ltimes \cC_{q}, \cC_{p})$. The
conjugate equals
\begin{equation*}
    \begin{bmatrix}
    1 & 0\\
    d^{-1} b a & 1\\
  \end{bmatrix}.
\end{equation*}
Since $d$ acts trivially on the quotient $(\cC_{p} \ltimes \cC_{q}) /
\cC_{q}$, we get that the automorphism group has structure
\begin{equation*}
  \cC_{p} \rtimes (\cC_{p-1} \times \Hol(\cC_{q})),
\end{equation*}
with $\cC_{p-1}$ acting as $\Aut(\cC_{p})$ and $\Hol(\cC_{q})$ acting
trivially, that is
\begin{equation*}
  \Hol(\cC_{p}) \times \Hol(\cC_{q}).
\end{equation*}

\bibliographystyle{amsalpha}
 
\providecommand{\bysame}{\leavevmode\hbox to3em{\hrulefill}\thinspace}
\providecommand{\MR}{\relax\ifhmode\unskip\space\fi MR }
\providecommand{\MRhref}[2]{%
  \href{http://www.ams.org/mathscinet-getitem?mr=#1}{#2}
}
\providecommand{\href}[2]{#2}

\end{document}